\documentclass{article}
\usepackage[utf8]{inputenc}

\usepackage{geometry}
\usepackage[utf8]{inputenc}
\usepackage{graphicx}
\usepackage{bbold}
\usepackage{amsmath}
\usepackage{amsthm}
\usepackage{amssymb}
\usepackage{float}
\usepackage{caption}
\usepackage{hyperref}
\usepackage{algorithm}
\usepackage{xcolor}
\usepackage{mathtools}
\usepackage{dsfont}
\usepackage{comment}
\usepackage[algo2e,ruled,vlined]{algorithm2e} 
\usepackage{algorithm}
\usepackage{algorithmicx}
\usepackage{stmaryrd}

\usepackage[noend]{algpseudocode}
\usepackage[backend=bibtex, style=alphabetic, maxbibnames=5]{biblatex}
\DeclareMathOperator{\R}{\mathbb{R}}

\def\Uc{{\cal U}}

\def \E{\mathbb{E}}

\def \Q{\mathcal{Q}}
\def \N{\mathbb{N}}
\def \P{\mathbb{P}}

\def\argmin{\mathop{\rm argmin}}

\newtheorem{Theorem}{Theorem}[section]

\newtheorem{Remark}[Theorem]{Remark}

\title{Reservoir optimization and Machine Learning methods}
\author{Xavier \textsc{Warin}
\footnote{EDF R\&D, FiME \sf \href{mailto:xavier.warin at  edf.fr}{xavier.warin at edf.fr}} 
}

\begin{document}
\maketitle
\begin{abstract}
Optimization of storage  using neural networks is now commonly achieved by solving a single optimization problem. We first show that this approach allows solving high-dimensional storage problems, but is limited by memory issues. We propose a modification of this algorithm based on the dynamic programming  principle and propose neural networks that outperform classical feedforward networks to approximate the Bellman values of the problem.
Finally, we study the stochastic linear case and show that Bellman values in storage problems can be accurately approximated using conditional cuts computed by a very recent neural network proposed by the author.
This new approximation method combines linear problem solving by a linear programming solver with a neural network approximation of the Bellman values.
\end{abstract}

\section{Introduction}
In industry, reservoirs are storage facilities  used to manage a stock of fuel in order to satisfy a particular objective function. 
For example, hydraulic reservoirs use water to generate electricity and the manager's goal is to supply energy to meet a demand at the lowest cost.  Gas storage is a special case of storage where the objective is to maximize profits by buying and selling gas on the market. 
Another example of storage is batteries, where electricity is directly injected or withdrawn, leading to systems that can be valued directly on the market or used to secure a global electrical system. \\
Traditional dynamic programming methods are generally used when the number of storage sites is small.
This dynamic programming method applied to a dynamical system breaks it down in time into a sequence of simpler problems. It provides, at a given time, the value of the system for optimal decisions taken from that time and from a {\bf state variable}, which is a set of parameters that fully describe the system.  This value is called the {\bf Bellman value} of the system and is therefore a function of time and other state variables.  In the case we face, the state of the system is fully described by a purely controlled deterministic state, the stock level in the storage, and the stochastic state (the {\bf uncertainties} in the sequel). This state of the system   is supposed to be {\bf Markovian}, meaning that what happens at a time $t+1$ depends only on the state at time $t$. \\
Gas storage  is often optimized and hedged using this type of approach \cite{warin2012gas}: the storage level is discretized on a grid, and Bellman values are calculated by regression using either the Bellman values at the following time step  \cite{tsitsiklis2001regression} or the cash generated at the following dates according to the classical Longstaff Schwartz approach \cite{longstaff2001valuing}.
These regressions are coupled with linear interpolation of the Bellman values or cash generated at different points on the grid. The method can only be used for non-linear problems in low dimension for two reasons. The most obvious one is the computational time, which explodes with the dimension, and computer clusters are needed to cope with this computational cost even with a number of reservoirs limited to 3 or 4. The second reason, which is in fact the first limiting one, is the need to store in memory the Bellman values needed by the software. This has led to the development of algorithms that split the  Bellman values into the memory of the different nodes of the computer cluster (see \cite{makassikis2007distribution} and the StOpt library \cite{gevret2018stochastic} for a recent implementation).  But even with this kind of approach, it is difficult to optimize problems with more than 4 or 5 reservoirs. Another pitfall that managers encounter even in dimension one, as it is the case for the optimization of gas storage facing the gas market, is the loss of concavity observed due to the regressions and interpolations in the calculation of conditional expectations, even when the solution is known to be concave with respect to the storage level.\\
In most countries with a large number of dams, dynamic programming methods are not used directly and the Stochastic Dual Dynamic Programming (SDDP) method  \cite{pereira1991multi} is generally used to manage the dams using a cut approximation of Bellman values, which are concave, for a linear objective function, with respect to the level of water in the reservoirs. Transition problems are solved using LP solvers with Bellman cuts as upper bounds on the final value. When uncertainties need to be incorporated into the state, breaking the concavity of the Bellman values, trees are often introduced and cuts need to be generated at each node of the tree, as explained in \cite{pereira1991multi}.  Another approach is to generate conditional cuts using regressions \cite{van2020conditional}. In all cases, forward passes (exploring the possible uncertainties and the levels of memory visited) and backward passes (adding cuts at the levels visited in the forward pass) give an iterative method that converges \cite{shapiro2011analysis}. An advantage of this method is that since cuts are used to generate an approximation to the Bellman values, concavity with respect to the storage level is preserved, and the marginal cost of the system (i.e., the derivative of the Bellman value with respect to the storage level) decreases with the storage level.
However, the method can be very slow to converge when the number of transition steps is large, and it can only be used for linear or special quadratic problems.\\
The use of neural networks to optimize  reservoir management  is old, as it was first used for gas storage in \cite{barrera2006numerical}, at a time when no public automatic differentiation software using neural networks was available, and the maximization of the profit generated by the storage had to be achieved using a gradient method that explicitly calculated the gradients. In this paper, feedforward networks are used to approximate the control (gas injection - withdrawal) at each time step and the solution obtained is compared with the solution obtained using a tree method. This part of the previous article was largely ignored by the scientific community at the time. More recently, this type of control representation has been used to solve BSDE problems \cite{han2018solving}, leading to much work on this approach.
Using the   dynamic programming principle, \cite{bachouch2021deep} proposed two algorithms using neural networks and gave some numerical results for the valuation of a storage :
\begin{itemize}
    \item The first, "Control Learning by Performance Iteration", has a Longstaff-Schwartz flavor: the control at the current time is approximated by a neural network, and the controls computed at the previous time iterations (so at the next time steps, since the process is backward) are reused to estimate the expectation of the objective function, starting from a randomized initial state at the current time. In practice, this method can only be used for a very small number of time steps N, as the cost of recalculating the objective function using the previously calculated controls is in $O(N^2)$. This approach has recently been successfully applied to deterministic control problems with very high order schemes to reduce the number of time steps used \cite{bokanowski2022neural}.
    On our test cases the method is not applicable.  
    \item The second method "Hybrid now" consists in computing the control at the current time step by a first optimization by neural networks and then estimating the Bellman values at the current time step by regression by a second optimization problem using a second feedforward network.
\end{itemize}
Finally, in a very recent article, \cite{curin2021deep} studies the valuation and hedging of a gas storage using a single optimization approximating the control at each time step using feedforward neural networks as in \cite{barrera2006numerical}. They also propose to "merge" the network between different time steps (i.e. to introduce a dependence on time in a network shared between different time steps). This kind of approach merging all time steps (so using a single neural network shared between all time steps) is the one proposed in \cite{fecamp2019risk} for risk evaluation or \cite{chan2019machine} for BSDE resolution: they show that it gives better stabilized results on these problems.\\
The present article first studies the "Global Valuation" method (GV) of one or more storages using the global approach used in \cite{barrera2006numerical} and \cite{curin2021deep} testing different formulations and network to represent the control.\\
This type of solution is impossible to use on real problems that involve optimization, for example, on the global year on an hourly basis for a storage, or on the global year even on a daily basis when many storages are interconnected: this is the case, for example, for the valuation  of a large number of batteries with a management that has an impact on the price of electricity.\\
In this case, we have to split the problem and use in a second section a modification of the "hybrid now" scheme, which solves the problem backwards, but estimates the control not on a single time step, but for a whole period. This scheme is called the "Global Split Dynamic Programming" (GSDP) method, and the "hybrid now" scheme is a special case where the transition problem is solved on a single time step. This scheme is studied on some difficult test cases and the use of near-optimal network to represent the Bellman values is studied.
\\
The GV method and the GSDP method are effective when non inter-temporal constraints have to be handled.
In the last section, we develop a methodology to handle intertemporal constraints in the special linear case: we develop a new algorithm GMCSDP to solve stochastic linear problems using a dynamic programming approach based on LP resolutions, as the SDDP method does, but with only a single backward pass. The Benders cuts approximating the Bellman values are generated using the new GroupMax neural network proposed in \cite{warin2022groupmax}.\\
The main results of this article are as follows:
\begin{itemize}
    \item 
    When the underlying process of the problem is Markovian, it is possible to optimize reservoirs very accurately with simple feedforward networks with a small number of layers and neurons for controls even in high dimension using the GV method. However, it is necessary to use a network for each time step and the "merged" network must be avoided as it leads to very poor results when the case is  stochastic. When the underlying process is not Markovian, we propose a combination of a LSTM network with classical feedforward networks to solve the problem. We show the effectiveness of the methodology.
    \item The use of the GSDP method results in a loss of optimality due to the calculation of Bellman values by regression with classical feed-forward networks. This loss of optimality is greatly increased when we use the "hybrid now" scheme. We show that with a single reservoir, it is possible to obtain reasonable results with a feedforward network by increasing the number of layers and neurons. Using 5 or 10 reservoirs optimized together, we show that the feedforward network is not effective for estimating Bellman values. To use the GSDP method in high dimension, we develop a network inspired by \cite{amos2017input} and show that it gives reasonable results in all dimensions tested, even with a small number of layers and neurons. When the problem is concave with respect to the storage level, we show that the network proposed in \cite{amos2017input} can be used to obtain good results while preserving concavity, and we show that this network is outperformed by an extension of the new GroupMax network proposed in \cite{warin2022groupmax}.
    \item   The GMCSDP method allows solving stochastic linear problems with inter-temporal constraints, but only in low and medium dimension.
\end{itemize}
In the whole article we will focus on a maximization of the profit of management in expectation.
This choice is driven by the fact that, in practice, it is the main concern associated with this type of management. The introduction of hedging strategies as proposed in \cite{curin2021deep} leads to the need to use a risk function to discriminate an optimal strategy. Another important point related to this choice is the fact that we can implement multi-storage optimization problems that can be easily reduced to a single-storage optimization, allowing to obtain a reference by the classical dynamic programming approach using regressions. Thus, we can verify that neural networks are indeed capable of solving problems of high dimension.

\section{The global approach for storage optimization}
First, for convenience, recall that a feedforward network with $K$ hidden layers and $m$ neurons is an operator $\phi$ $:$
$\R^{d_0} \longrightarrow \R^{d_1}$
defined by the following recurrence:
\begin{align}
z_0 = & x  \in \R^{d_0},\\
z_{i+1}= & \rho( \sigma_i z_i+ b_i),  0 \le i  <  K, \\
\phi(x)= & \hat \rho( \sigma_{K} z_{K} + b_{K}), 
 \label{eq:feedforward}
\end{align}
where :
\begin{itemize}
    \item $\sigma_1 \in  \R^{m \times d_0}$, $\sigma_i \in  \R^{m \times m}$ for $i=1, \dots,K-1$, $\sigma_K \in  \R^{d_1 \times m}$, 
    \item $b_i \in \R^m$,  for $i=0, \dots,K-1$, $b_K \in \R^{d_1}$,
    \item $\rho$, $\hat \rho$ are non linear activation functions (Elu, Relu, tanh etc..) applied component wise.
\end{itemize}
The set $\theta$ of the network parameters is defined by the set of all matrix and vector coefficients of $\sigma_i$, $b_i$, $i=0, \dots, K$.\\
First, we test different approximations and formulations for a linear storage problem.
Retaining the best formulation,
we then extend our tests to higher dimensions by considering a nonlinear problem that combines the management of different reservoirs.
\subsection{On a linear problem in dimension one}
Suppose we manage a storage (gas storage, battery) on a commodity market (gas, electricity) where the commodity follows the HJM model
\begin{align}
    \frac{dF(t,T)}{F(t,T)} = e^{-a(T-t)} \sigma d W_t
    \label{eq:price}
\end{align}
where $W_t$ is a one-dimensional Brownian motion defined on a probability space $(\Omega, \mathcal{F}, \P)$. The spot price is then $S_t= F(t,t)$.\
The characteristics of the storage are the withdrawal $C_W$ and injection $C_I$ rates (both positive) during a time step $\Delta t$, its maximum capacity $Q_{Max}$, and its initial inventory $Q_{Init}$.\\
The reservoir manager wants to maximize the expected profit associated with filling the reservoir by buying the commodity when prices are low and emptying the reservoir to sell the commodity when prices are high.
We then define the objective function for $N$ optimization dates $t_i = i \Delta t$, for $i=0, \dots, N-1$:
\begin{align}
    J(U)= -\E[ \sum_{i=0}^{N-1} S_{t_i} u_i],
    \label{eq:J1}
\end{align}
where  $ U= (u_i)_{i=0,N-1}$ is  in the set $\Uc$ of the   non anticipative admissible strategies such that:
\begin{align}
    0 \le Q_j := Q_{init} +  \sum_{i=0}^{j-1} u_i \le Q_{Max}, \quad  \mbox{ for } j=0,\dots, N ,  \nonumber\\
     -C_W \le u_i \le C_I , \quad  \mbox{ for } j=0, \dots,N-1.
     \label{eq:flow}
\end{align}
We  want to maximize the expected gain associated with storage management:
\begin{align}
    J^{*} = \sup_{U \in \Uc} J(U).
\end{align}
Let us define $Q$ as the reservoir level.
There are many ways to deal with the constraints imposed on the level of the storage ($Q$ must remain positive and below $Q_{max}$): among them, clipping the control, penalizing the objective function are possible, but the best approach (we won't report results on less effective approaches) consists in using the \cite{curin2021deep} approach.
First we introduce for a given $i$ in $0,\dots, N-1$:
\begin{align}
     \hat C_I^i= ((Q_i+ C_I) \wedge Q_{max}) -Q_i,    & \quad  \hat C_W^i=  Q_i - ((Q_i- C_W) \vee 0).
     \label{eq:const}
\end{align}
Since the problem is Markov in $(S,Q)$, we can introduce as in \cite{barrera2006numerical} a feed-forward network $\phi_i^{\theta_i}$ with parameters $\theta_i$ per time step $i$ as an operator from $\R^2$ to $[0,1]$ (using a sigmoid activation function at the output and a $\tanh$ activation function for hidden layers)  such that the control is approximated by:
\begin{align*}
     -\hat C_W^i + (\hat C_W^i +\hat C_I^i) \phi_i^{\theta_i}.
\end{align*}

Noting  $\theta = (\theta_i)_{i=0,N-1}$, we approximate the optimal reservoir management by solving :
\begin{align}
    \theta^*= \argmin_{\theta} \E[ \sum_{i=0}^{N-1} S_{t_i} (-\hat C_W^i + (\hat C_W^i +\hat C_I^i) \phi_i^{\theta_i}(S_{t_i},Q_i) ) ],
    \label{eq:optNNN}
\end{align}
where $Q_i$ follows \eqref{eq:flow}, $\hat C_W^i$ and $\hat C_I^i$ are given by \eqref{eq:const}, and the dynamics of $F$ follows \eqref{eq:price}.\\

A second version classically consists of introducing a single neural network $\phi$ ("merged" network) with parameters $\theta$ as a function of $(t,F,Q)$ and modifying the optimization \eqref{eq:optNNN} as follows:
\begin{align}
    \theta^*= \argmin_{\theta} \E[ \sum_{i=0}^{N-1} S_{t_i} \big(-\hat C_W^i + (\hat C_W^i +\hat C_I^i) \phi^{\theta}(t_i,S_{t_i},Q_i) \big) ].
    \label{eq:optOneNN}
\end{align}

We use a classical stochastic gradient descent ADAM method in Tensorflow \cite{2015tensorflow} and use normalized data for $F$ and $Q$ as input to the neural network. \\
For the test case, we assume that we are optimizing a storage for $N=365$ days with one decision per day ($\Delta t =1$). The price parameters (expressed in days) are $\sigma=0.08$, $a=0.01$. The initial forward curve has seasonal and weekly variations, as in the energy market, and is given by
$F(0,T)= 30 +5 \cos(\frac{2 \pi T}{N}) + \cos(\frac{2 \pi T}{7})$.
For the reservoir we take $C_W=10$, $C_I=5$, $Q_{Max}=100$, $Q_{Init}=50$.
 A reference is computed using the StOpt library \cite{gevret2018stochastic} by dynamic programming with adaptive linear regression \cite{bouchard2012monte} and cash flow interpolations as exposed in \cite{warin2012gas}. The optimization uses $100$ basis functions and $10^6$ trajectories to optimize the regression control. The parameters are taken such that the problem is hard to solve by dynamic programming using regressions. Since the solution is bang bang \cite{barrera2006numerical}, we use $20$ grid points to discretize the storage and only bang bang controls are tested leading to a very fast estimation. Then a simulation is performed using the Bellman values obtained in the optimization. The value obtained in the optimization is equal to $4938$, while the value in the simulation, taken as a reference, is equal to $4932$.
\begin{figure}[H]
    \centering
    \begin{minipage}[c]{.49\linewidth}
    \includegraphics[width=\linewidth]{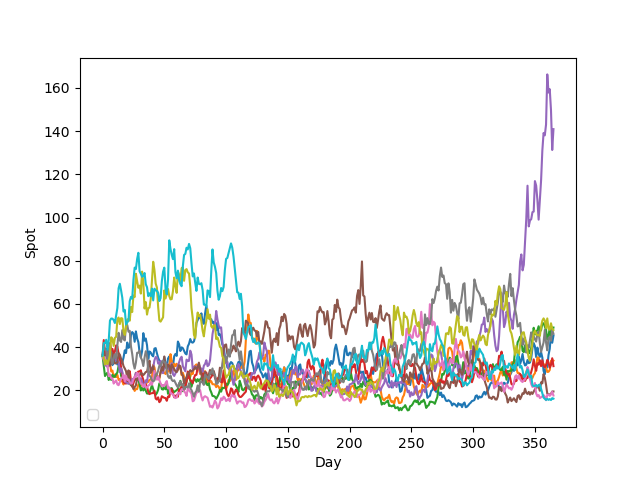}
    \caption*{Prices}
    \end{minipage}
    \begin{minipage}[c]{.49\linewidth}
    \includegraphics[width=\linewidth]{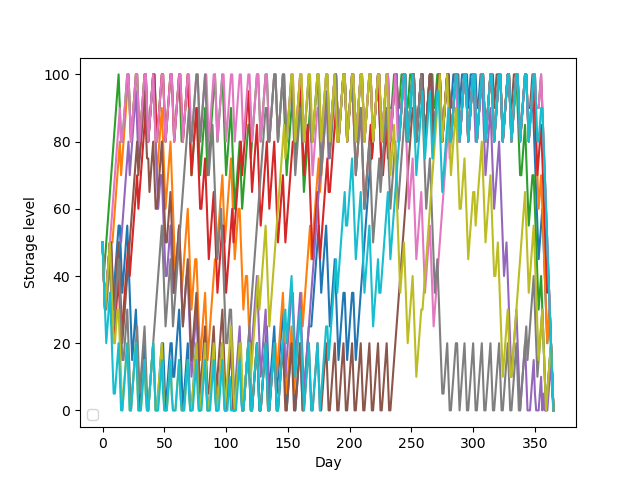}
    \caption*{Optimal storage trajectories}
    \end{minipage}
    \caption{10 spot and  optimal management trajectories.  }
    \label{fig:LionSDeriv500}
\end{figure}
Using 2 hidden layers with 11 neurons, we train the problems \eqref{eq:optNNN} and \eqref{eq:optOneNN} using a minibatch of size 200, 100000 iterations for the gradient descent, and an initial learning rate equal to $2 \times 10^{-3}$.
The results obtained after training with 200000 trajectories are given in the table \ref{tab:global1D} using 10 runs.
Minimal/Maximal is the minimal/maximal value obtained on the 10 runs, while Average is the average value obtained.
\begin{table}[H]
    \centering
    \begin{tabular}{|c|c|c|c|c|}
        \hline
  &  Maximal & Minimal & Average &  Min diff with DP \\ \hline
 One network per day \eqref{eq:optNNN} & 4925  & 4914 & 4922 &  7 \\ \hline
 A singe network  \eqref{eq:optOneNN} & 3944  &1795 & 3702 & 988 \\ \hline
    \end{tabular}
    \caption{ Neural network valuation with 10 runs. }
    \label{tab:global1D}
\end{table}
Results are excellent with a  network by day but very poor with a single network. The results with a single network don't change while increasing the number of layers and neurons.
Using a network per time step, we can check that the results remain very good in all cases while decreasing the volatility $\sigma$.
\subsection{Linear problem increasing the dimension}
\label{sec:linnD}
To get a linear problem in a higher dimension, we assume that we have $M$ similar storages to manage, each with a strategy $(u^j_0, \dots, u_{N-1}^j)$ for $j=1, \dots, M$. 
We note  $Q_i= (Q_i^j)_{j=1,\dots,M}$, where $Q^j_i$ is the level in storage $j$ at time $t_i$, $U= ((u^j_i)_{i=0,N-1})_{j=1,M}$ and  the function to maximize is given by:
\begin{align}
    J^M(U)= -\E[ \sum_{j=1}^M \sum_{i=0}^{N-1} S_{t_i} u_i^j],
    \label{eq:linEq}
\end{align}
and
\begin{align}
    J^{M,*} = \sup_{U \in \Uc} J^M(U),
    \label{eq:maxProb}
\end{align}
where all strategies satisfy a flow equation similar to \eqref{eq:feedforward}.\\

Similarly to the previous section, we introduce a neural network per time step depending on the current prices and the different storage levels. The network   with parameters $\theta_i$ per time step $i$, as an operator from $\R^{1+M}$ to $\R^M$,  approximates  the controls $(u_i^1, \dots, u_i^M)$ using the same activation functions as in the previous section. It leads to optimize:
\begin{align}
    \theta^*= \argmin_{\theta}  \E[ \sum_{i=0}^{N-1} S_{t_i} (-\hat C_W^i + (\hat C_W^i +\hat C_I^i) \phi_i^{\theta_i}(S_{t_i},Q_i)). \mathbb{1}_M  ],
    \label{eq:optNNNB}
\end{align}
where now the $\hat C_W^i$ and $\hat C_I^i$ are now vectors in $\R^M$ with each component satisfying an equation similar to \eqref{eq:const}.\\

In this section, we test two networks:
\begin{itemize}
    \item First, the classical feedforward network $\phi$ previously introduced,
    \item Secondly, since the solution is symmetric, we test the DeepSet network \cite{deepsets} (with the same parameters for the number of layers and neurons  as originally proposed by the authors) allowing to impose that the control satisfies the symmetry:
    $$ u^l(S, Q^1, \dots , Q^l, \dots,  Q^m,\dots, Q^M) = u^m( S, Q^1, \dots , Q^m, \dots,  Q^l,\dots, Q^M),$$ for all $(Q^1,\dots, Q^M)$ state positions in the storages.
    This network has proven to be more effective than feedforward networks, but mainly for very high-dimensional PDEs arising from a special approximation of mean-field problems \cite{germain2021deepsets}. 
\end{itemize}
\begin{Remark}
Since the storage have the same characteristics,  the global value of the set of storages is the same if, for example, storages $i$ and $j$ permute their levels $Q^i$ and $Q^j$. This invariance by permutation of the storage level allows us to search for symmetric control.
\end{Remark}
We further assume that, at the initial date, all storages have the same state: then solving \eqref{eq:maxProb} is equivalent to solving the problem for one storage and we have $J^{M,*}= M J^{1,*}$.
For each storage, we take the same storage characteristics as before.
The training  and simulation  parameters are the same as in the previous section except that the number of neurons used is  $10+M$.
The table \ref{tab:globalMD} shows the results obtained by the two networks. The results obtained by the classical feedforward networks are already optimal, and the DeepSet networks are not interesting in such small dimensions.
\begin{table}[H]
    \centering
    \begin{tabular}{|c|c|c|c|c|c|}
        \hline
  Network & M &  Maximal & Minimal & Average &  Min diff with DP \\ \hline
  feedforward & 3 & 4926  & 4915 & 4921 &  6 \\ \hline
  feedforward &10 & 4931  & 4918 &  4925 &  1 \\ \hline
  DeepSet &  3 & 4896 & 4882 &  4889 & 35 \\ \hline
  DeepSet & 10 &  4904&  4891 & 4896 & 28  \\ \hline
    \end{tabular}
    \caption{ Neural network valuation divided by the dimension in the linear case. 10 runs.}
    \label{tab:globalMD}
\end{table}
\subsection{Results on a  non linear case}
\label{sec:nonLin}
To get reference results for a non-linear case, we now assume that the price is no longer exogenous and that the impact is proportional to $\frac{P}{M}$, which leads to a modification of the function to maximize:
\begin{align}
    J^M(U)= -\E[ \sum_{j=1}^M \sum_{i=0}^{N-1} (S_{t_i}+ \frac{P}{M} \sum_{l=1}^M u_i^l) u_i^j].
\end{align}
This kind of modeling of a price impact is necessary, for example, in battery management, when the amount of energy managed represents a fairly large fraction of the energy available in the market.\\ 
Using the same methodology as before, taking $P=0.2$ and the same parameters as in the previous sections, we can first for $M=1$ obtain a reference with a very thin discretization of the command and the grid storage using dynamic programming with the StOpt library. The value obtained by classical regression is equal to $3802$ in optimization and a value in simulation, taken as a reference, is equal to $3796$. 

The solution of \eqref{eq:maxProb} again satisfies  $J^{M,*}= M J^{1,*}$ and we give  $\frac{J^{M,*}}{M}$ obtained by the previously defined feed forward network for this non linear case  in the table \ref{tab:globalMDNLinear}.
\begin{table}[H]
    \centering
    \begin{tabular}{|c|c|c|c|c|}
        \hline
  M &  Maximal & Minimal & Average &  Min diff with DP \\ \hline
  1 &  3794 & 3780  & 3788 &  7 \\ \hline
  5 &  3799   & 3789 &  3794 &  2\\ \hline
  10 &  3797 & 3784 &  3791 &  5\\ \hline
    \end{tabular}
    \caption{ Neural network valuation with ten run  divided by the dimension in the non linear case.}
    \label{tab:globalMDNLinear}
\end{table}
Again, the results are very good in dimension 1 to 10 for this very stochastic case, as shown in the table \ref{tab:globalMDNLinear}: the error with respect to the reference computed by dynamic programming is nearly zero and hardly distinguishable from the Monte Carlo error associated with the reference. Therefore, the GV method provides results comparable in accuracy to very efficient dynamic programming methods. 

\subsection{Extension in the  non Markovian case}
\begin{figure}[H]
    \centering
    \includegraphics[width=6cm]{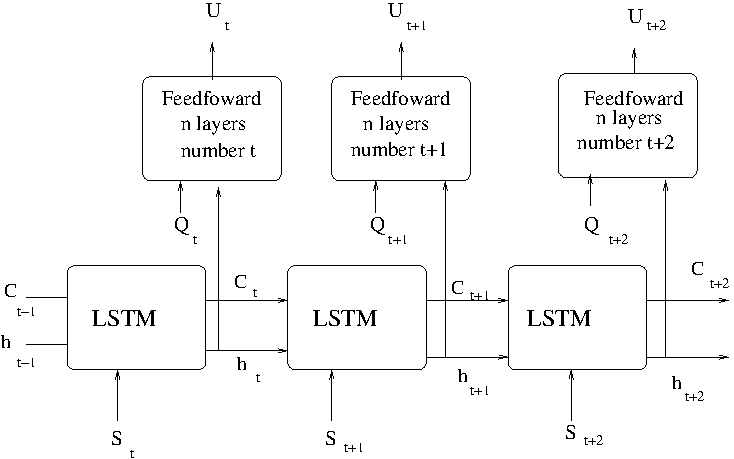}
    \caption{Unrolled LSTM with feedforward to approximate control.}
    \label{fig:LSTM}
\end{figure}

If the price, or more generally the uncertainties, are not Markovian, it is possible to extend the previous feedforward network to deal with this feature.\\
Since the optimal control, at a time $t$, is a function of the entire history of the price $(S_u)_{u \le t}$ and the current position in the storage, the idea is to use a recurrent network such as a LSTM network \cite{hochreiter1997long} to deal with the price dependence. At each time step, the output of the LSTM network (with the price as input) and the position in the storage are used as inputs to a feedforward network that provides the control.
The figure \ref{fig:LSTM} shows an unrolled version of the LSTM connected to the feedforward at each time step.
 To test this network we now assume that the future price \eqref{eq:price} is  replaced by \eqref{eq:pric3D}.
\begin{align}
    \frac{dF(t,T)}{F(t,T)} =  \sum_{i=1}^3 e^{-a_i(T-t)} \sigma_i d W_t^i,
    \label{eq:pric3D}
\end{align}
where now $W_t =(W^1_t, W^2_t, W^3_t)$ is a three dimensional Brownian motion.\\
We still consider the linear problem \eqref{eq:J1} with the same characteristics as before, but using the parameters $(\sigma_1, \sigma_2, \sigma_3)= (0.04, 0.028, 0.023)$ and  $(a_1, a_2, a_3)= (0.01,0.005,0.0033).$
In the table \ref{tab:globalMDLSTM}, we compare the results obtained with the feedforward that takes as input the three risk factors of the price model $( \sigma_i\int_0^t e^{-a_i (t-s)} dW^i_s)_{i=1,3}$ and the stock level with the results obtained with the LSTM feedforward, where the LSTM network takes as input the price history and has as output $50$ units. Using dynamic programming with the StOpt library, we were able to obtain a value of $4300$, achieving an optimization regression with the \cite{bouchard2012monte} method in dimension 3 with $10^7$ trajectories and $10^3$ meshes.  We were unable to refine the results due to memory constraints.

\begin{table}[H]
    \centering
    \begin{tabular}{|c|c|c|c|c|}
        \hline
  Network & M &  Maximal & Minimal & Average  \\ \hline
  feedforward & 1  &  4332 & 4322  &  4328  \\ \hline
  feedforward & 5 &  4333 & 4320 & 4327 \\ \hline
    feedforward &  10 &  4334 & 4318 & 4329\\ \hline
  LSTM-feedforward  & 1 & 4285 & 4277 &  4280 \\ \hline
  LSTM-feedforward  & 5 &  4284& 4273 & 4279   \\ \hline
  LSTM-feedforward  & 10 &4285  & 4268 &  4279 \\ \hline
    \end{tabular}
    \caption{ Neural network valuation divided by the dimension in the linear case for the non Markovian case (10 runs).}
    \label{tab:globalMDLSTM}
\end{table}
\begin{Remark}
The LSTM method is not needed in this simple case because we can increase the state by taking 3 future values with different maturities to recover a Markov state.
\end{Remark}
The results are very good with a very small loss of accuracy compared to the feedforward network. We note that the results obtained with the feedforward networks are slightly better than the results obtained with our not converged classical regression method. Here, the feedforward neural network methodology provides even better results than classical dynamic programming methods. The proposed LSTM feedforward network also gives accurate results, but slightly below the dynamic programming results using regression. 

\section{Global Split Dynamic Programming method}
\label{sec:GSDP}
The global method proposed in the previous section is very effective but may be impossible to implement when the number of dates is too large: memory problems appear  and another approach has to be used. This leads to the development of a combination of the Hybrid-Now method of \cite{bachouch2021deep} and the global method of \cite{barrera2006numerical}.
Suppose that the objective function is: 
\begin{align}
    J(U) = \E[ \sum_{i=0}^{N-1} f(t_i, S_{t_i}, U_i)],
    \label{eq:genOpProib}
\end{align}
where the control $U_i$ is a vector of size $M$,
with flow constraints as in \eqref{eq:flow} but applied to the control component by component. We assume that $S_{t}$ is Markov.\\
Suppose that we want to solve \eqref{eq:maxProb} and 
that we divide the $N$ dates in $N= \sum_{l=1}^L \hat N_l$ with $\hat N_l \in \mathds{N}^{*}$ for $l=1, \dots, L$.\\
The idea is to use the dynamic programming principle to replace the previous global optimization by $L$ backward optimizations where each optimization $l$ is achieved  on $\hat N_l$ time steps.
Then an optimization  on the period $[t_{(\sum_{j=1}^{l-1} \hat N_j)-1}, t_{(\sum_{j=1}^{l} \hat N_j)-1}]$  allows us to obtain the Bellman values for  each possible state at time $t_{(\sum_{j=1}^{l-1} \hat N_j)-1}$. This Bellman value is then used if $l>1$ as the final value of the optimization on $[t_{(\sum_{j=1}^{l-2} \hat N_j)-1}, t_{(\sum_{j=1}^{l-1} \hat N_j)-1}]$.
\\
The algorithm \ref{algo:GSDP} allows to solve the problem \eqref{eq:maxProb}.\\
\begin{algorithm2e}[H]
 	\caption{GSDP method }
 	\label{algo:GSDP}
 	\textbf{Output:} Estimates the Bellman values at dates $t_0$,  $t_{\sum_{l= 1}^{\tilde l} \hat N_l}$, for $\tilde l =1, L-1$, and all optimal controls.\\
 	
 	$\tilde i = N$\\
 	$VB^L= 0$\\
 	\For{$l$ $=$ $L,\ldots,1$}{
 	    $\tilde i = \tilde i - \hat{N}_l$\\
  		Introduce $\hat N_l$ feedforward networks  $\phi_k^{\theta_k}$ on $\R^{M+1}$ with values in $\R^M$  with sigmoid activation at the output,
  	    \begin{flalign}
        \theta^*:= (\theta_0^*,\dots, \theta_{\hat N_l-1}^*)= \argmin_{\theta}  \E[ \sum_{i=0}^{\hat N_l-1} f(t_{\tilde i+i},S_{t_{\tilde i +i}}, U_i^\theta) + VB^{l}( S_{\tilde i + \hat{ N}_l},Q_{\hat N_l})], 
        \label{eq:contrHyb}
        \end{flalign}
        such that $Q_0 \sim U[0, Q_{max}]^M$  and for $0 \le i < \hat N_l$:
        \begin{flalign}
        \hat C_I^i= ((Q_i+ C_I) \wedge Q_{max}) -Q_i,&  \quad \quad    \hat C_W^i=  Q_i - ((Q_i- C_W) \vee 0)  \nonumber\\
        U_i^\theta= -\hat C_W^i + (\hat C_W^i +\hat C_I^i) \phi_i^{\theta_i}(S_{t_{\tilde i+i}},Q_i), & \quad
        Q_{i+1} =  Q_{i}+    U_i^\theta.
        \label{eq:flowAl}
        \end{flalign}  
        Introduce a neural network $\psi^\kappa$ with identity output activation function with parameters $\kappa$,
        \begin{align}
            \kappa^* = \argmin_{\kappa} \E[ \big( \sum_{i=0}^{\hat N_l-1} f(t_{\tilde i+i},S_{t_{\tilde i +i}}, U_i^{\theta^*})  + VB^{l}( S_{\tilde i + \hat{ N}_l},Q_{\hat N_l}) - \psi^{\kappa}(S_{t_{\tilde i}}, Q_0) \big)^2],
            \label{eq:regHyb}
        \end{align}
        where \eqref{eq:flowAl} is satisfied and $Q_0 \sim U[0, Q_{max}]^M$. \\
        $VB^{l-1} = \psi^{\kappa^*} $
     }
 \end{algorithm2e}
 \begin{Remark}
 In this  method it is necessary to obtain a representation of the Bellman values at each date $t_{(\sum_{j=1}^{l-1} \hat N_j)-1}$.
In the proposed  algorithm these Bellman values are estimated by introducing another neural network that is optimized by solving the equation \eqref{eq:regHyb}.
 \end{Remark}
 In this algorithm the use of a feedforward to solve  \eqref{eq:regHyb} may seem natural,  as it was proposed for the hydrid-now method in \cite{bachouch2021deep}.
 We take our  linear test case  in dimension 1 and test the use of this network for $\psi^\kappa$ to solve \eqref{eq:regHyb} for different number of $L$ values (only  $\hat N_0$ can be  different from the $\hat N_l$ for $l>0$, which are all equal :$\hat N_l =\hat N_m$ for $l>0$ and $m>0$).
 As for the neural networks used for the control in \eqref{eq:contrHyb}, we keep the same characteristics as in the previous sections.\\
 Using different numbers of neurons and layers to approximate the Bellman values, in the table \ref{tab:linearFeedGSDP} we give the results obtained with the best of  3 runs with $100000$ gradient iterations with an initial learning rate equal to $5 \times 10^{-3}$, still using the ADAM method. The activation function $\rho$ in \eqref{eq:feedforward} is a Relu function, since it gives better results than the $\tanh$ activation function, and it gives similar results to the $Elu$ activation function. 
 
 \begin{table}[H]
    \centering
    \begin{tabular}{|c|c|c|c|c|}
        \hline
  L & m & $\tilde L$ & solution &  Min diff with DP \\ \hline
  4  &   11  & 2 &  4900 &  32\\ \hline
  13  &  11 &  2 &   4816&  116\\ \hline
  53 &  11 & 2 &  4389  &   543\\ \hline
4  &   30  & 3 &  4899 &  33 \\ \hline
  13  &  30 &  3 &   4834&  98 \\ \hline
  53 &  30 & 3 &  4633  &  299 \\ \hline
    \end{tabular}
    \caption{ GSDP method for the one dimensional linear case using a feedforward  network using a number of neurons $m$, and a number of layers $\tilde L$ to solve \eqref{eq:regHyb}. We take the best result out of 10 runs.}
    \label{tab:linearFeedGSDP}
\end{table}
There is a loss of optimality increasing as $L$ increases. For $L=53$ we need to use at least 3 layers and 30 neurons to avoid too much loss of accuracy and keep acceptable results.
Using more layers or neurons can  improve the results very slightly and now we take $5$ layers with $20$ neurons for this  feedforward network. We  now test the previous  linear case and the non linear one from sections  \ref{sec:linnD} and \ref{sec:nonLin} keeping $L=53$ and letting the dimension increase. We report the  results in the table \ref{tab:dimFeedNeural}.
\begin{table}[H]
    \centering
    \begin{tabular}{|c|c|c|c|c|c|} \hline
      Test case   &  M &  max & min & average & min error \\ \hline
       Linear  & 1  &  4651 &  4543 &  4606&   280 \\ \hline
       Linear &  5  &  4352 & 4252 &  4179 & 579  \\ \hline
       Linear & 10 &   2547  & 1831       & 2151        &  1248     \\ \hline
       Non linear &  1 & 3687 & 3646  & 3664 &  108 \\ \hline
       Nonlinear  & 5 & 3074 & 2122 & 2529 & 721 \\ \hline
       Nonlinear & 10 & 2039 & 1379 & 1165 & 1756 \\ \hline
    \end{tabular}
    \caption{Feedforward  network results $\frac{J^{M,*}}{M}$ with the GSDP method for linear and non linear cases  with $5$ layers, $20$ neurons for \eqref{eq:regHyb}, $L=53$ and varying  the dimension $M$. We use $10$ runs and report the best (max), the worst (min) results, and finally the error for the best result obtained.}
    \label{tab:dimFeedNeural}
\end{table}
The deterioration of the result obtained with the dimension is obvious and even in dimension 5 the error obtained is far too important.
\begin{Remark}
The Hybrid-Now method is the limit case when $L=365$:  since the results are already bad for $L=53$, we can say that the method is ineffective to optimize real storage for a whole year.
\end{Remark}
The treatment of storage levels must be made different from the treatment of uncertainty, and it leads to the development of more adapted networks.
We  propose to use three different networks.
\subsection{A first network $\psi^A$ preserving concavity}
Since in this case, the solution is concave with respect to the stock level in the storage we can use a modification of the \cite{amos2017input} network  given by the  recursion \eqref{eq:modAmos}. Suppose that  the input to the neural network is $\tilde x= (x,y) \in \R_{d_0}$ where we have concavity in $y \in \R^{k}$,
\begin{flalign}
u_{i+1} = &  \tilde \rho( \tilde W_i u_i + \tilde b_i),  \nonumber \\ 
z_{i+1} = & \rho( [W^{(z)}_i \otimes (W^{(zu)}_i u_i + b^{(z)}_i) ]^+  z_i   + \nonumber \\ & W^{(y)}_i( y \circ (W^{(yu)}_i u_i + b^{(y)}_i)) + W^{(u)}_i u_i + b_i), \quad  \mbox{ for } i \le  K \nonumber \\
\psi^A = &z_{K+1},  \quad  u_0 =x , \quad z_0=0,
\label{eq:modAmos}
\end{flalign}
where  $\rho$ is a concave non increasing activation function which we take to be equal to  minus Relu, $\circ$ denotes the Hadamard product, $\otimes$ is applied between a matrix $A \in \R^{m \times n}$ and a vector $B \in \R^n$ such that $A \otimes B \in \R^{m \times n}$ and $(A \otimes B)_{i,j}=  A_{i,j} B_j$.  As noted  in \cite{amos2017input}, the concavity of the solution is given by the properties of $\rho$ and the fact that the weight before $z_i$ is positive in \eqref{eq:modAmos}.
Using $m_x$ neurons for the non-concave part of the function and $m_y$ neural networks for the convex part of the network, $\tilde W_0 \in \R^{m_x \times d_0-k}$, $\tilde W_i \in \R^{m_x \times m_x}$ for $i>0$,  $W^{(zu)}_i \in \R^{m_y \times m_x}$ for $i > 0$, $W^{(z)}_i \in \R^{m_y \times m_y}$ for $j<K$, $W^{(z)}_K \in \R^{1 \times  m_y}$ as the output is a scalar function. We don't go into detail about the size of the different matrices $W^{(y)}$, $W^{(yu)}$, $W^{(u)}$  and the different biases that are obvious. 

In all experiments, $\tilde \rho$ is the ReLU activation function.

\subsection{A second network  $\psi^{AD}$ removing the concavity constraints}
In a natural way, we modify the previous neural network by removing the concavity constraints  that lead to:
\begin{flalign}
u_{i+1} = &  \tilde \rho( \tilde W_i u_i + \tilde b_i),  \nonumber \\ 
z_{i+1} = & \rho( W^{(z)}_i  (z_i \circ (W^{(zu)}_i u_i + b^{(z)}_i))  +  \nonumber\\ & W^{(y)}_i( y \circ (W^{(yu)}_i u_i + b^{(y)}_i)) + W^{(u)}_i u_i + b_i), \quad  \mbox{ for } i \le  K,\nonumber \\
\psi^{AD} = & z_{K+1}, \quad   u_0 =x, \quad z_0=0.
\label{eq:modAmosNoConc}
\end{flalign}
 This gives us   a  neural network that can be used  for non concave/convex problems.
 In all experiments $\tilde \rho$ and $\rho$ are ReLU activation functions.
 
\subsection{The GroupMax  network $\psi^{GM}$  using cuts when the solution is concave}
The GroupMax is a recently developed neural network \cite{warin2022groupmax} that combines the ideas in \cite{anil2019sorting}, \cite{tanielian2021approximating}, and those in \cite{amos2017input}, but allows cuts to represent a concave solution. If the function is concave only with respect to $y$, the following neural network is suggested in \cite{warin2022groupmax} and  generates cuts conditional on $x$:
\begin{flalign}
u_0 =& x ,  \quad z_0= 0 ,\nonumber \\
u_{i+1} = &  \tilde \rho( \tilde W_i u_i + \tilde b_i),  \nonumber \\ 
z_{i+1} = & \rho( [W^{(z)}_i \otimes (W^{(zu)}_i u_i + b^{(z)}_i) ]^+  z_i   +  \nonumber\\ & W^{(y)}_i( y \circ (W^{(yu)}_i u_i + b^{(y)}_i)) + W^{(u)}_i u_i + b_i), \quad  \mbox{ for } i \le  K-1, \nonumber \\
\psi^{GM}(x,y)= &  \hat\rho( [W^{(z)}_K \otimes (W^{(zu)}_K u_K + b^{(z)}_K) ]^+  z_K   + \nonumber \\ & W^{(y)}_K( y \circ (W^{(yu)}_K u_K + b^{(y)}_K)) + W^{(u)}_K u_K + b_K),
\label{eq:groupMax}
\end{flalign}
where all the matrices involved have the same size as the matrices in \eqref{eq:modAmosNoConc} except that the matrix $W^{(z)}_K$ is in  $\R^{m_y \times  m_y}$.\\
In \eqref{eq:groupMax}, the $\tilde \rho$ is a classical activation function like ReLU and  in order to get conditional cuts to approximate the solution, we take $ \hat \rho$ as an activation function working on the whole vector:
\begin{align*}
    \hat \rho(x) = \min_{i=1, \dots, d} x_i \quad \mbox{for } x \in \R^d.
\end{align*}
The $\rho$ is defined by grouping the elements of the vector as in the GroupSort neural network \cite{anil2019sorting}. Assuming  $x \in \R^{m_y}$, $G  \le m_y$  $\in \N^{*}$  the group size such that $ \tilde m = \frac{m_y}{G} \in \N^{*}$  is  the number of groups, $\rho$ maps $\R^{m_y}$ to $\R^{\tilde m}$ such that:
\begin{flalign*}
\rho(x)_i = \min_{ j =1, \dots,G} x_{(i-1) G +j},\quad  \mbox{ for } i =1, \dots, \tilde m. 
\end{flalign*}
In \cite{warin2022groupmax}, it is shown that this network generates conditional cuts with respect to $x$.
In all experiments, $\tilde \rho$ is a ReLU activation function.\\
\subsection{Numerical results}
We test the three networks on the linear and the non linear case. It is obvious that in the linear case, the Bellman value is concave with  respect to the storage level. In the non linear case, the concavity  is still present \cite{girardeau2015convergence}. 
All the Bellman values obtained by the three networks are not very sensitive to the number of layers and neurons. For the three networks used to approximate the Bellman values, we take $m_x=10$ and  $3$ hidden layers, taking $m_y=20$ for the first two  networks and $m_y=40$ for the GroupMax network. We keep the same parameters as in the previous sections to estimate the controls.
\begin{table}[H]
    \centering
     \centering
    \begin{tabular}{|c|c|c|c|c|c|} \hline
      M & Network   &   max & min & average & min error \\ \hline
  1 &  $\phi^A$&    4679  &  4618 & 4645 & 252 \\ \hline
 1 & $\phi^{AD}$ & 4798 & 4738 & 4771 & 133 \\ \hline
 1  & $\phi^{GM}$ & 4810&  4777 & 4795 &  122 \\ \hline
  5 &  $\phi^A$&   4352   & 3949 & 4179 & 579 \\ \hline
 5 & $\phi^{AD}$ &  4482& 4318 & 4399  & 449  \\ \hline
 5  & $\phi^{GM}$ & 4641& 4519  & 4601  & 290  \\ \hline
   10 &  $\phi^A$& 4027   & 3724   & 4027 & 904 \\ \hline
 10 & $\phi^{AD}$ & 4151 & 3890 & 4031 & 780 \\ \hline
 10  & $\phi^{GM}$ & 4425&  4316 & 4351 & 506  \\ \hline
     \end{tabular}
    \caption{Result $\frac{J^{M,*}}{M}$ of the Linear case with the GSDP method with L=53 for the different networks. 10 runs.}
    \label{tab:GSDPNewNetLin}
\end{table}

\begin{table}[H]
    \centering
     \centering
    \begin{tabular}{|c|c|c|c|c|c|} \hline
      M & Network   &   max & min & average & min error \\ \hline
  1 &  $\phi^A$&    3585  &  3540  & 3558 & 210  \\ \hline
 1 & $\phi^{AD}$ & 3687 & 3646 & 3664 & 108 \\ \hline
 1  & $\phi^{GM}$ & 3688& 3663  & 3676 &  107 \\ \hline
  5 &  $\phi^A$&   3444   & 3206 & 3412 & 351  \\ \hline
 5 & $\phi^{AD}$ &  3538   & 3278 & 3434 & 257    \\ \hline
 5  & $\phi^{GM}$ &  3633  & 3589 & 3614  &  162  \\ \hline
   10 &  $\phi^A$& 3456   & 3044    & 3288 & 339 \\ \hline
 10 & $\phi^{AD}$ & 3389 & 3308 & 3205 &  406\\ \hline
 10  & $\phi^{GM}$ & 3556&  3482 & 3556 & 239  \\ \hline
     \end{tabular}
    \caption{Result $\frac{J^{M,*}}{M}$ on the Non Linear case with the GSDP method with L=53 for the different networks. 10 runs.}
    \label{tab:GSNNewNetNonLin}
\end{table}

As the dimension $M$ increases, the variance of the results obtained increases. The GroupMax is clearly superior to the other networks.
The results are better, the loss of accuracy decreases more slowly with the dimension and the variance of the results obtained is much lower  than with the other networks.\\
We conclude  that it is optimal to use a number $L$ as small as possible. When the problem is not concave with respect to the storage levels, the $\phi^{AD}$ network is the best, and when the problem is concave, the cut methodology given by the GroupMax network $\phi^{GM}$ is the best choice.\\
We have shown that we are able to  efficiently optimize   joint reservoirs  even with a large number of stocks and  in the non Markovian case.
To circumvent the memory limitation  of the GV method, the GSDP method can be used but the neural network used to approximate the Bellman values has to be chosen carefully while  keeping  the number of global sub-resolutions  as small as possible.
Depending on the concavity of the Bellman values, a neural network may be preferred over the other. All these methods are effective  as long as the position in the storage is not constrained to be discrete and as long as there are no inter-temporal constraints.\\
When it is crucial to take  these inter-temporal constraints into account, other methods can be proposed  in the  simpler linear case, as we develop in the next section.

\section{GroupMax Cut Split Dynamic Programming  (GMCSDP) method}

High-dimensional reservoir management is often accomplished using a linear stochastic model.
Using a hazard decision framework, solving the problem is based on  a stochastic model where uncertainties are revealed for a period $[t,t+\Delta t[$ and decisions are made on each of the $P$ sub-intervals $\delta t$ of $\Delta t$.
This type of model is often used to optimize, for example, energy producing assets with inter-temporal constraints on the period $[t,t+\Delta t[$. For example, uncertainties are to be revealed every week (time step $\Delta t$). Then the dynamic programming approach is used with a weekly time step, and the optimal orders during the week for each hour (time step $\delta t$) are computed either using a deterministic dynamic programming approach (if the coupling between assets is not too difficult to account for) or using a linear programming solver if the problem is linear. \\
If the problem is linear, using the fact that the Bellman values are concave with respect to the storage levels, on each interval $[t_j,t_{j+1}]$ where $t_j= j \Delta t$, for $j=1, \ldots,N$, injection/withdrawal decisions are made at each date $t_j + k \delta t$, $k=0,\ldots, P-1$, starting from a level in the reservoirs $Q=(Q_1, \ldots, Q_M)$ at date $t_j$. This optimization is solved using an LP problem with end conditions given by cuts of the Bellman values (called Benders cuts) at date $t_{j+1}$ (see for example the recent article by \cite{leclere2020exact} for details in the SDDP approach).
Since the Bellman values are generally not concave with respect to the uncertainties, the cuts are conditional on the uncertainty level. These cuts are often given at the nodes of a scenario tree \cite{pereira1991multi} and can also be computed by regressions as in \cite{van2020conditional}.    \\
\begin{itemize}
    \item When the dimension is low, starting point in the reservoirs for the LP problems are derived by exploring a grid, and a single backward pass is achieved. This is a dynamic programming method.
    \item As the dimension increases, an iterative process of backward pass and forward exploration simulations is used. The forward pass allows the reservoir levels of interest to be revealed and then exploited during the subsequent backward pass, which adds new cuts. This iterative procedure is called SDDP \cite{pereira1991multi}.
\end{itemize}
The two methods are developed in open source  libraries such as StOpt \cite{gevret2018stochastic} using regressions and trees.
These methods are very popular because they allow to handle difficult constraints. These constraints are taken into account by the LP solver used to solve the transition problems.\\
\begin{Remark}
The nonlinear transition problem with a quadratic concave cost function with respect to the control and the reservoir level can also be solved by a quadratic solver using the same methodology, but at a higher cost.
\end{Remark}
As mentioned in the introduction, convergence can be very slow and the stopping criterion be can be difficult to implement, especially when the Bellman values are not concave with respect to the uncertainties.\\
The use of dynamic programming methods to estimate the cuts on a lattice suffers from both computational time and memory requirements. We propose to use the GroupMax network to estimate the Bellman values by cuts.\\
Suppose that we want to solve equation the following equation
\begin{align}
J^* = & \sup_{U=(U_0,U_{N-1}) \in \Uc} \sum_{i=0}^{N-1} \E[ f(t_i, U_i,S_{t_i}) ],  \nonumber \\
\underline U \le &  U_i \le \bar U, \mbox{ with } U_i \in \R^p ,  \quad i=0, N-1,\nonumber\\
X_0 = &\tilde X \in \R^q,  \nonumber\\
X_{i+1} =& X_i + A_i(S_i) U_i + B_i(S_i) \in \R^q, \quad i=0, N-1, \nonumber\\
\underline X \le &  X_i \le \bar X, \quad i=0, N,
\label{eq:genLP}
\end{align}
where $f$ is linear or quadratic concave with respect to $U$, $S_t$ is a discrete time Markov process in $\R^d$, $A_i$ a function from $\R^d$ to $ \R^{q \times p}$, $B_i$ function from $\R^d$ in  $\R^q$ for $i=0,N-1$.\\
In a water storage optimization problem, $S_i$ would represent stochastic inflows and the electricity price at sub-interval $[t_i, t_{i+1}[$, while $X_i$ would represent the storage level at each date of $[t_i, t_{i+1}[$ depending on the turbine commands $U_i$ at $[t_i, t_{i+1}[$ and the inflows giving the flow equation involving $X_i$ above. $f$ would be, as before, a linear function of the command $U_i$ where the coefficients are a function of the electricity price as in equation \eqref{eq:J1}.\\
Based on the dynamic programming principle and using the GroupMax network, we propose to use the algorithm \ref{algo:GMCSDP} to optimize \eqref{eq:genLP}, where the Bellman values are estimated by cuts.

\begin{algorithm2e}[H]
 	\caption{GMCSDP method }
 	\label{algo:GMCSDP}
 	\textbf{Output:} \\
 	
 	$VB^{N-1}= 0$\\
 	\For{$i$ $=$ $N-1,1$}{
 	    Introduced a GroupMax neural network $\psi^{GM, \theta}(S,X)$ with parameter $\theta$
 	    \begin{align*}
 	        \theta^* = \argmin_{\theta} \E[ (\psi^{GM, \theta}(S_{t_{i-1}},X) - \Q(S_{t_i},X))^2],
 	    \end{align*}
 	    where
 	    \begin{align}
 	    \label{eq:sddpLike}
 	        \Q(S_{t_{i}},X)= &\max_{U \in \R^p} f(t_i,S_{t_i},U) +  \xi, \nonumber \\
 	        \xi \le &  VB^{i}(S_{t_i},\tilde X),\\
 	        \tilde X = & X + A_i(S_{t_i}) U + B_i(S_{t_i}), \nonumber\\
 	        \underline X \le &  \tilde X \le \bar X,  \nonumber\\
 	        \underline U \le &  U \le \bar U, \nonumber
 	    \end{align}
 	    with $ X \sim U([\underline X, \bar X])$, $S_{t_{i-1}}$ sampled from $S_0$ and $S_{t_i}$ sampled from $S_{t_{i-1}}$. \\
 	    $VB^{i-1}(s,x) = \psi^{GM,\theta^*}(s,x)$
     }
     Optimize first time step:
     \begin{align*}
 	        &\max_{U \in \R^p} f(0,S_0,U) +  \xi, \\
 	        \xi \le &  VB^{0}(S_{0},\tilde X),\\
 	        \tilde X = & X_0 + A_0(S_{0}) U + B_0(S_{0}),\\
 	        \underline X \le &  \tilde X \le \bar X, \\
 	        \underline U \le &  U \le \bar U.
 	    \end{align*}
 \end{algorithm2e}
 To test the algorithm, we suppose that the problem is linear, given by \eqref{eq:maxProb}, \eqref{eq:linEq}
 with the flow equation \eqref{eq:flow}.
 The characteristics of the storages are unchanged. The initial forward curve is given by $F(0,T)= 30 + 4 \cos(\frac{2 \pi T}{7})$. The price model is still given by \eqref{eq:price} but with the parameters $\sigma=0.3$, $a=0.16$. 
 We take $N=42$ and the optimization by the dynamic programming method using the property that the optimal control are bang bang gives a value of $3426$ while the value obtained in forward using the optimal control is $3424$.
 \begin{Remark}
 In this test, at each time step, the transition problem involves only one subtime step, and one can easily test the commands to find the optimal one and obtain the optimal asset value with a simple dynamic programming algorithm using the StOpt library and the Longstaff Schwarz method.
 \end{Remark}
 Using the algorithm \ref{algo:GMCSDP}, in a first test case we use a GroupMax network with 2 layers (one hidden layer) and a group size equal to 2. 
 The ADAM stochastic gradient descent is used with a batch size of $200$, a number of gradient iterations of $15000$, and an initial learning rate of $5 \times 10^{-3}$, decreasing to $10^{-4}$. Local optimizations are done with the Coin LP solver.
 Results for $\frac{J^*}{M}$ are given in the table \ref{tab:GMCSDP}.

 \begin{table}[H]
     \centering
     \begin{tabular}{|c|c|c|c|c|c|} \hline
         M & $m_y$   &   max & min & average & min error   \\ \hline    
         1 & 8  &  3362  & 3317 &  3335 &    61 \\ \hline
        1 & 10  & 3370  &  3321 &  3349 &  53   \\ \hline
         1 & 12  & 3355  &  3334 &  3345  &   69  \\ \hline
         3   & 8 & 3207 & 3113 & 3163 & 218  \\ \hline
         3  &  10 &   3225    &   3173    &   3202   &  200   \\ \hline
        3  &  12  &   3233    & 3166      & 3194     &  192    \\ \hline
          5 & 8  &  3063 &  2758 &  2913 &   360  \\ \hline
         5 & 10  &  3012 &  2729 & 2885  &   411    \\ \hline
         5 & 12  &  3052 &  2799   &   2980 &    371      \\ \hline
    \end{tabular}
     \caption{GMCSDP results on 10 runs for $\frac{J^*}{M}$. Reference is SDP with regressions.}
     \label{tab:GMCSDP}
 \end{table}
 The accuracy decreases as the dimension increases, while the variance of the result obtained increases.
 We don't see any clear differences using different values of $m_y$, except perhaps in dimension 3, where increasing $m_y$ seems to give slightly better results.
 \begin{Remark}
 Since the time induced by the LP solving is strongly related to the number of cuts used, we limit the number of layers to 2 ($K=1$, which allows to have $m_y 2^{\frac{m_y}{2}}$ cuts using $m_y$ neurons).  We keep $m_x$ equal to $8$.
 As shown in \cite{warin2022groupmax}, it is necessary to increase the number of layers to get high accuracy, but this leads to very time-consuming problems to solve.
 This increase in the computational cost of the LP is due to the fact that the equation \eqref{eq:sddpLike} adds a constraint for each cut of the terminal Bellman value. A reduction of the computational cost would be possible by pruning (eliminating) the inactive constraints \cite{pfeiffer2012two}, which are numerous (see graphs of cuts in \cite{warin2022groupmax}).
 \end{Remark}
 At last we take a smaller test case,  using a storage  with $C_I=10 \times  \mathbb{1}_M $, $C_W=20 \times  \mathbb{1}_M $, $Q_{Max}=100 \times  \mathbb{1}_M $, $Q_{init}=50  \times \mathbb{1}_M $.
 We keep the same parameters  $\sigma=0.3$, $a=0.16$ for the price model. The forward curve is $F(0,T)= 30 + 4 \cos(\frac{2 \pi T}{4})$ and $N=8$.
 We test the influence  of the number of neurons $m_y$, the group size $G$, and  test a  number of  layers $K$ equal to 2 or even 3 in low dimension. We take  $m_x=6$.  
  The results are given in the table \ref{tab:GMCSDP2}. The reference value is $1818$.
 \begin{table}[H]
     \centering
     \begin{tabular}{|c|c|c|c|c|c|} \hline
         M & $K$  & $m_y$   & $G$ &   Solution & Error   \\ \hline    
         1 & 2 & 9 & 3  & 1803 &  15\\ \hline
         1 & 2 & 10 & 2  & 1807 &  11\\ \hline
         1 & 2 & 12 & 2  & 1809 &  9\\ \hline
         1 & 3 & 10 & 5  & 1790 &  28\\ \hline
         2 & 2 & 9 & 3  & 1795 &  23\\ \hline
         2 & 2 & 10 & 2  & 1797 &  21\\ \hline
         2 & 2 & 12 & 2  & 1794 &  24\\ \hline
         2 & 3 & 10 & 5 &  1780 &  38 \\ \hline
         3 & 2 & 9 & 3  & 1753 &  65\\ \hline
         3 & 2 & 10 & 2  & 1767 &  51\\ \hline
         3 & 2 & 12 & 2  & 1774 &  44\\ \hline
         4 & 2 & 9 & 3  & 1728 &  89\\ \hline
         4 & 2 & 10 & 2  & 1736 &  82\\ \hline
         4 & 2 & 12 & 2  & 1727 &  91\\ \hline
    \end{tabular}
     \caption{GMCSDP best results on 10 runs for $\frac{J^*}{M}$. Reference is SDP with regressions.}
     \label{tab:GMCSDP2}
 \end{table}
 We observe the same behavior as in the previous case: the results are very accurate in dimension one and there is a degradation as we increase the dimension.
 We test three layers in the special case where $m_y=10$, $G=5$ in low dimension, which significantly increases the computation time, but with results slightly below those obtained with two layers. Changing the group size and the number of neurons with three layers remains too costly even in this simple test case.
 Using a group size of 3 instead of 2 with 2 layers doesn't seem to change the results much.
 Therefore, the proposed approach seems to be limited to low or medium dimensions, and in order to use it efficiently, two layers with 10 or 12 neurons and a group size of 2 seems to be optimal. \\

  In real industrial problems, constraints between storage reduce the volatility of the system and convergence should be  easier to achieve.
 This approach allows to avoid the memory cost due to the storage of the Bellman functions and only the computing time remains a constraint: the elimination of inactive cuts, the parallelization by threads and MPI should allow to efficiently reduce this computational cost.

\section{Conclusion}
We have proposed an effective method to optimize reservoirs with new neural networks even in high dimension without inter-temporal constraints, allowing to circumvent the limitation of the Global Value method. We have also proposed a method that could be used to optimize linear problems taking into account the inter-temporal constraints in a multistage stochastic framework: the results indicate that a very high number of cuts must be used to approximate Bellman values, limiting the potential of the method. Pruning methods could be used to reduce the time taken by the linear programming solvers, then allowing more than 2 layers in the GroupMax network with different group size and number of neurons. The method can be used to solve low or medium dimensional problems.
Nevertheless, the nonlinear case with inter temporal constraints remains a challenge for neural networks: classical penalties, usually added to the objective functions to impose these difficult constraints, tend to give bad results by smoothing the solution.\\
Another open question is the possibility of solving this type of problem with neural networks when the stock in the reservoir is constrained to take only discrete values.

\printbibliography 

\end{document}